\documentclass[journal]{IEEEtran}
\usepackage{lineno,hyperref}
\usepackage{amsfonts}
\usepackage{amssymb}
\usepackage{amsmath}
\usepackage{graphicx} 
\usepackage{enumerate}
\usepackage{booktabs} 
\usepackage{hypernat} 
\usepackage{stfloats}

\newcommand \T {\mathrm{T}}

\usepackage{cite}

\ifCLASSINFOpdf
\else
\fi
\begin{document}
%
\title{Contraction Analysis on Primal-Dual Gradient Optimization}
%
%
%

\author{Yanxu~Su,
	~Yang~Shi,~\IEEEmembership{~Fellow,~IEEE,}
	and~Changyin~Sun
	\thanks{This work was supported in part by the National Natural Science Foundation of China under Grant 61520106009 and Grant U1713209.}
	\thanks{Y. Su and C. Sun are with the School of Automation, Southeast University, Nanjing 210096, China, and also with the Key Laboratory of Measurement and Control of Complex System of Engineering, Ministry of Education,
		Southeast University, Nanjing 210096, China (e-mail: yanxu.su@seu.edu.cn; cysun@seu.edu.cn).}
	\thanks{Y. Shi is with the Department of Mechanical Engineering, University of Victoria, Victoria, BC V8W 3P6, Canada (e-mail: yshi@uvic.ca).}}

%
%

\markboth{Journal of \LaTeX\ Class Files,~Vol.~14, No.~8, August~2015}%
{Shell \MakeLowercase{\textit{et al.}}: Bare Demo of IEEEtran.cls for IEEE Journals}
%



\maketitle

\begin{abstract}
This paper analyzes the contraction of the primal-dual gradient optimization via contraction theory in the context of discrete-time updating dynamics. The contraction theory based on Riemannian manifolds is first established for convergence analysis of a convex optimization algorithm. The equality and inequality constrained optimization cases are studied, respectively. Under some reasonable assumptions, we construct the Riemannian metric to characterize a contraction region. It is shown that if the step-sizes of the updating dynamics are properly designed, the convergence rates for both cases can be obtained according to the contraction region in which the convergence can be guaranteed. Moreover, the augmented Lagrangian function which is projection free is adopted to tackle the inequality constraints. Some numerical experiments are simulated to demonstrate the effectiveness of the presented contraction analysis results on primal-dual gradient optimization algorithm. 
\end{abstract}


%
\IEEEpeerreviewmaketitle

\section{Introduction}
%
%
%
%
\IEEEPARstart{T}{he} following constrained optimization problem is studied in this paper
\begin{subequations}\label{OP}
	\begin{align}
	\min_{x\in\mathbb{R}^n}\quad &f\left(x\right)\\
	\text{s.t.} \quad &A_1x=b_1\\
	\quad &A_2x\le b_2
	\end{align}
\end{subequations}
where the objective function $f\left(x\right)$ is strongly convex and smooth with respect to $x\in\mathbb{R}^n$, and the constraints are depicted by $A_1\in\mathbb{R}^{p_1\times n}$, $A_2\in\mathbb{R}^{p_2\times n}$, $b_1\in\mathbb{R}^{p_1}$ and $b_1\in\mathbb{R}^{p_1}$. Given the Lagrangian multiplier $\lambda\in\mathbb{R}^{p_i}$ for $i=1,2$, the Lagrangian (or Augmented Lagrangian) associated with optimization problem in (\ref{OP}) is $\mathcal{L}\left(x,\lambda \right)$. The Primal-Dual Gradient Optimization (PDGO) with discrete-time dynamics under study in this paper is given as follows
\begin{subequations}\label{PDGD}
	\begin{align}
	x^{k+1} &= x^k -\alpha \nabla_{x} \mathcal{L}\left(x^k,\lambda^k \right) \\
	\lambda^{k+1} &=\lambda^k + \beta \nabla_{\lambda} \mathcal{L}\left(x^k,\lambda^k \right)
	\end{align}
\end{subequations}
where $k$ is the iteration number and the positive scalars $\alpha$ and $\beta$ are step-sizes.

\subsection*{Background}

Convex optimization is of great interest in recent years for its wide implementations in control systems ranging from distributed optimization \cite{shi2015extra,xi2017dextra}, power grid \cite{yi2015distributed,li2015connecting}, machine learning \cite{mairal2015incremental}, game theory \cite{gharesifard2013distributed} and so on. Various convex optimization algorithms have been developed to provide better performance, such as first-order algorithms \cite{nedic2009distributed,shi2015extra}, dual gradient algorithms \cite{patrinos2013accelerated,wu2019fenchel}, Alternating Direction Method of Multipliers (ADMM) algorithms \cite{shi2014linear,majzoobi2019analysis}, primal-dual gradient algorithms \cite{hale2017asynchronous,bernstein2019online}, etc. The convergence analysis of an optimization algorithm plays a crucial role in evaluating the performance. The asymptotic convergence of the primal-dual gradient optimization has dawn tremendous attentions in the existing literature \cite{cherukuri2016asymptotic,cherukuri2017saddle,goebel2017stability}. In practice, however, the exponential convergence (also called geometric convergence) is a desired property for its stronger convergence guarantees. The geometric convergence of primal-dual gradient with continuous-time dynamics has been investigated \cite{nguyen2018contraction,qu2018exponential,liang2019exponential}, while it has been rarely studied in the context of discrete-time dynamics. In particular, the PDGO with discrete-time dynamics in (\ref{PDGD}) is essential when it is adopted as a computational tool in the convex optimization.

The contraction theory is a widely adopted approach for nonlinear systems analysis. In accordance with the fluid mechanics and differential geometry, the contraction theory was first established in \cite{Lohmiller1998}. The traditional approach based on Riemannian manifolds has been generalized to many fields such as distributed nonlinear systems \cite{long2018distributed}, stochastic incremental systems \cite{pham2009contraction}, etc. In addition, some recent results were published inspired by Finsler manifolds \cite{forni2014differential,chaffey2018control}. It, however, is worth to mention that the existing results were obtained for continuous-time dynamical systems. To the best of our knowledge, the contraction theory has been rarely studied in the context of discrete-time dynamical systems. Moreover, only a few results have been addressed using the contraction theory to analyze the convergence of an optimization algorithm \cite{nguyen2018contraction}. The objective of this paper is to properly design the step-sizes via contraction theory towards the geometric convergence guarantees of PDGO in (\ref{OP}).

\subsection*{Literature review}

In \cite{qu2018exponential}, the authors investigated the exponential convergence of primal-dual gradient optimization in the framework of continuous-time dynamics. Furthermore, they extended their theoretical analysis results to the discrete-time dynamics by using Euler discretization. In \cite{liang2019exponential}, the authors studied the exponential convergence of the continuous-time distributed primal-dual gradient optimization without strong convexity. They generalized their analysis results to the discrete-time dynamics via Euler’s approximation method. It is worth to note that the convergence analysis under the discrete-time dynamics relies on the theoretical analysis results for the continuous-time dynamics, which may lead to conservative.

In \cite{cortes2018distributed}, a class of saddle-point-like dynamics was studied and a novel Lyapunov function was constructed to demonstrate the exponential convergence rate guarantees of the primal-dual optimization algorithms subject to equality constraints. However, the theoretical results cannot be generalized to inequality constrained optimization directly.

In \cite{nguyen2018contraction}, the strict contraction of continuous-time primal-dual gradient optimization was analyzed by means of contraction theory. Moreover, the robustness of the PDGO was exploited in specific metrics. The discrete-time dynamics, however, were not considered.

\subsection*{Contribution}

In this paper, we investigate the constrained optimization problem in (\ref{OP}) with primal-dual gradient updating dynamics in (\ref{PDGD}). Under some reasonable assumptions, the geometric convergence of the discrete-time PDGO in (\ref{PDGD}) is rigorously analyzed by using contraction theory. Meanwhile, the convergence rates of equality and inequality constrained optimization cases are given, respectively. The proofs rely on the contraction analysis for discrete-time dynamics based on Riemannian manifolds. Notice that the contraction theory which is a well-known method for nonlinear system analysis is first adopted to analyze the geometric convergence of an optimization algorithm. The Riemannian metric is constructed to depict the contraction region in which the geometric convergence can be guaranteed. Furthermore, a classic augmented Lagrangian function which is projection free is adopted to handle the inequality constrained optimization problem. 

\subsection*{Outline}

The remainder of this paper is structured as follows. Section \ref{Preliminaries} gives some crucial preliminaries for convergence analysis of an optimization algorithm. The theoretical analysis results of the equality and inequality constrained optimization cases are summarized in Section \ref{Main Results}. In Section \ref{Numerical Examples}, some numerical experiments are simulated to verify the effectiveness of the proposed primal-dual gradient optimization algorithm. Section \ref{Conclusion} concludes the paper.

\section{Preliminaries}\label{Preliminaries}

\subsection{Notational Conventions}

The notations adopted in this paper are stated in the following. Denote real and natural number set as $\mathbb{R}$ and $\mathbb{N}$, respectively.  A matrix with the superscript $\T$ represents its transposition. The superscript $+$ of a state vector stands for the successor state. The subscripts of $\mathbb{N}_{\left[a,\,b\right]}$ is the integers in the interval $\left[a,\,b\right]$. We denote $\left[z\right]_+ = \max\left\lbrace z,0 \right\rbrace$. The column vector with the $i$-th element be $1$ is represented as $\text{e}_i$. The identity matrix $\boldsymbol{I}$ in this paper is with proper dimension. 

\subsection{Contraction Theory}

We recap the crucial definitions of Riemannian geometry in the following. For more details, please refer to \cite{Lohmiller1998,liu2017robust} and references therein. For two vectors $\delta_1, \delta_2$ on the tangent space of a given manifold, the Riemannian metric is a smoothly varying inner product $\left\langle \delta_1, \delta_2 \right\rangle_x = \delta_1^\T M\left(x\right) \delta_2$ with respect to a positive matrix function $M\left(x\right)$. Throughout this paper, we use the notation as $\|\delta\|_{x} = \sqrt{\left\langle \delta_1, \delta_2 \right\rangle_x}$. Notice that the matrix $M\left(x\right)=M$ is assumed to be constant in this paper. Given a pair of points $x\in\delta_1$ and $y\in \delta_2$, let $\Gamma\left(x, y\right)$ be a smooth path connecting $x$ and $y$, which implies that there exists a piecewise smooth mapping $\varrho\in\Gamma\left(x, y\right)$ satisfying $\varrho\left(0\right) = x$ and $\varrho\left(1\right) = y$. The Riemannian length is defined as $L\left(\varrho\right):=\int_{0}^{1} \|\varrho_s\|_{\varrho} ds$, the Riemannian energy $E\left(\varrho\right):=\int_{0}^{1} \|\varrho_s\|_{\varrho}^2 ds$, where $\varrho_s := \partial\varrho_s/\partial s$. Denote the Riemannian distance as $d\left(x,y\right) = \inf_{\varrho\in\Gamma\left(x,y\right)}L\left(\varrho\right)$. Without loss of generality, we define $E\left(x,y\right) := d\left(x,y\right)^2$ in this paper. 

Some vital results on contraction theory for convergence analysis of PDGO with discrete-time dynamics are briefly stated in the following. For continuous-time dynamical systems, the readers are referred to \cite{Lohmiller1998} for detailed introduction of contraction analysis. Consider an autonomous discrete-time nonlinear dynamical system described by
\begin{equation}\label{autonomousSys}
\xi\left(t+1\right)=\varPhi\left(\xi\left(t\right),t\right),
\end{equation}
where $\xi\left(t\right)\in \mathbb{R}^{n}$ represents the state vector at time instant $t\in \mathbb{N}$ and $\varPhi$ is a smooth and differentiable function. Denote the differential dynamics of the system in (\ref{autonomousSys}) as 
\begin{equation}\label{diffDynamic}
\delta \xi\left(t+1\right) = \frac{\partial \varPhi}{\partial \xi}\delta \xi\left(t\right).
\end{equation}
Let the target state trajectory be a forward-complete solution of the system in (\ref{autonomousSys}). If there exists a controller such that 
\begin{equation}\label{locExpStable}
d\left(\xi\left(t\right),\xi^{\star}\left(t\right)\right) \le C \lambda ^{t}d\left(\xi\left(0\right),\xi^{\star}\left(0\right)\right)
\end{equation}
for $t\in\mathbb{N}$, where $\lambda$ is the convergence rate and $C$ is a positive scalar which are independent of the initial states, the target state trajectory $x^{\star}$ is said to be globally exponentially controllable. 

In accordance with the contraction region defined in \cite{Lohmiller1998}, we can similarly establish the following lemma for discrete-time nonlinear systems.

\textit{Lemma 1:} Given a discrete-time nonlinear system $\xi\left(t+1\right)=\varPhi\left(\xi\left(t\right),t\right)$, a region in state space is called a contraction region, if there exists a uniformly positive definite constant metric $M$, such that
\begin{equation}\label{contractionTheoryDef}
\frac{\partial \varPhi^{\T}}{\partial \xi}M\frac{\partial \varPhi}{\partial \xi} - M \le \left(\tau^2-1\right) M
\end{equation}
with the convergence rate $\tau\in\left(0,1\right)$.

\textit{Proof:} For the given discrete-time nonlinear system $\xi\left(t+1\right)=\varPhi\left(\xi\left(t\right),t\right)$, the differential dynamics can be expressed as (\ref{diffDynamic}). Suppose that there exist two state trajectories with different initial values, which can be referred to the actual $\xi$ and the target $\xi^{\star}$. Denote the tangent vector of the actual state trajectory as $\delta \xi$. By introducing the constant Riemannian metric $M$, we can obtain
\begin{equation}\label{contractionTheoryDefProof1}
\begin{aligned}
&\left(\delta \xi^+\right)^{\T}M\delta \xi^+-\delta \xi^{\T}M\delta \xi\\
=&\delta \xi^{\T}\frac{\partial \varPhi}{\partial \xi}^{\T}M\frac{\partial \varPhi}{\partial \varPhi}\delta \xi-\delta \xi^{\T}M\delta \xi\\
=&\delta \xi^{\T}\left(\frac{\partial \varPhi}{\partial \xi}^{\T}M\frac{\partial \varPhi}{\partial \xi}-M\right)\delta \xi\\
\le&\delta \xi^{\T}\left(\tau^2-1\right)M\delta \xi
\end{aligned}
\end{equation}
where $\tau$ is the maximum eigenvalue of $\frac{\partial \varPhi}{\partial \xi}$. According to (\ref{locExpStable}), it can be guaranteed that the geometric convergence of the Riemannian distance within the contraction region depicted by (\ref{contractionTheoryDef}) with $\tau \in \left(0,1\right)$. The proof is completed. $\hfill\blacksquare$

\subsection{Optimization}

Some rational assumptions are given as follows, which are adopted in this paper to analyze the convergence of the constrained optimization problem in (\ref{OP}).

\textit{Assumption 1:} The objective function $f\left(x\right)$ under study in (\ref{OP}) is twice differentiable and $\underline{\rho}$-strongly convex with Lipschitz gradient, such that
\begin{equation}\label{Assuption1}
\underline{\rho} \boldsymbol{I} \le \left \langle \nabla f\left(x_1\right)-\nabla f\left(x_1\right),x_1-x_2 \right \rangle \le\overline{\rho} \boldsymbol{I}
\end{equation}
for all $x_1,x_2\in\mathbb{R}^n$, where $\overline{\rho}$ is the Lipschitz constant satisfying $0<\underline{\rho}\le\overline{\rho}<\infty$.

\textit{Assumption 2:} Given the matrices $A_1$ and $A_2$ with full row rank, there exist $0<\underline{\sigma}_1\le\overline{\sigma}_1<\infty$ and $0<\underline{\sigma}_2\le\overline{\sigma}_2<\infty$ such that 
\begin{subequations}\label{Assuption2}
	\begin{align}
	\underline{\sigma}_1\boldsymbol{I}\le A_1^\T A_1 \le\overline{\sigma}_1\boldsymbol{I}\\
	\underline{\sigma}_2\boldsymbol{I}\le A_2^\T A_2 \le\overline{\sigma}_2\boldsymbol{I}
	\end{align}
\end{subequations}

\textit{Remark 1:} Assumption 2 is qualified by the linear independence property of the constraints, which is standard in convergence analysis of an optimization algorithm \cite{qu2018exponential}. Furthermore, Assumptions 1 and 2 guarantee the solution (saddle-point) of primal-dual gradient optimization to be unique.

\section{Main Results}\label{Main Results}

In this section, the equality and inequality constrained optimization cases are studied separately. The Riemannian metric is constructed to characterize the contraction region within which the geometric convergence is guaranteed. Moreover, the convergence rates of both cases are given according to the contraction theory.

\subsection{Equality Constrained Optimization}

Consider the equality constrained optimization problem as follows
\begin{subequations}\label{ECO_OP}
\begin{align}
	\min_{x\in\mathbb{R}^n}\quad &f\left(x\right)\\
	\text{s.t.} \quad &A_1 x=b_1 
\end{align}
\end{subequations}
where $f$ and $A_1$ satisfy Assumptions 1 and 2. The Lagrangian of the considered optimization problem in (\ref{ECO_OP}) is given in the following
\begin{equation}\label{LangrangianECO}
\mathcal{L}_{EC}\left(x,\lambda\right) = f\left(x\right)+\lambda^\T\left(A_1 x-b_1 \right)
\end{equation}
where $\lambda\in\mathbb{R}^{p_1}$ is the Lagrangian multiplier. The Karush-Kuhn-Tucker (KKT) conditions associated to the problem in (\ref{ECO_OP}) for characterizing the optimal pair $\left(x^*,\lambda^* \right)$ can be described by
\begin{subequations}\label{KKT}
	\begin{align}
		\nabla_{x^*} \mathcal{L}_{EC}\left(x^*,\lambda^*\right) &= \nabla f\left(x^*\right)+A^\T \lambda^* = 0\\
		\nabla_{\lambda^*} \mathcal{L}_{EC}\left(x^*,\lambda^*\right) &= A_1 x^*-b_1 = 0
	\end{align}
\end{subequations}
Thus, the primal-dual gradient dynamics depicting the optimization updating can be expressed by
\begin{subequations}\label{PDGD_ECO}
	\begin{align}
	x^{k+1} &= x^k -\alpha \nabla_{x} \mathcal{L}_{EC}\left(x^k,\lambda^k\right)\nonumber \\
	&= x^k-\alpha\nabla f\left(x^k\right)-\alpha A_1^\T \lambda^k\\
	\lambda^{k+1} &=\lambda^k + \beta \nabla_{\lambda} \mathcal{L}_{EC}\left(x^k,\lambda^k\right)\nonumber\\
	& =\lambda^k+ \beta \left(A_1 x^k-b\right)
	\end{align}
\end{subequations}
where $\alpha$ is the step-size for primal updating and $ \beta $ is for dual updating. The following lemma is introduced to analyze the convergence of the primal-dual gradient optimization problem in (\ref{ECO_OP}).

\textit{Lemma 2:} \cite{qu2018exponential} Under Assumption 1, given the optimal state $x^*$, it holds that 
\begin{equation}\label{G_bound}
f\left(x\right)-f\left(x^*\right) = G\left(x\right)\left(x-x^*\right)
\end{equation}
where $G\left(x\right)$ is a symmetric matrix with respect to $x$, satisfying $\underline{\rho}\le G\left(x\right)\le\overline{\rho}$ for any $x\in\mathbb{R}^n$.

In what follows, the convergence of the primal-dual gradient dynamics is rigorously analyzed and the theoretical analysis results are summarized in the following theorem.

\textit{Theorem 1:} Suppose Assumptions 1 and 2 hold. For the primal-dual gradient optimization updating dynamics in (\ref{PDGD_ECO}), if it holds that $\max\left\lbrace\frac{\beta\overline{\sigma}_1}{\underline{\rho}},\alpha\overline{\rho} \right\rbrace \in\left[0,\frac{1}{2}\right]$ for properly designed step-sizes $\alpha$ and $\beta$, there exist two constants $C_1$ and $C_2$ which depend on $\underline{\rho}$, $\overline{\rho}$, $\underline{\sigma}_1$, $\overline{\sigma}_1$, $\alpha$ and $\beta$, such that the following conditions
\begin{subequations}\label{ECO_Theorem}
\begin{align}
	d\left(x^k, x^* \right) \le C_1 \tau_{EC}^k d\left(x^0, x^* \right)\\
	d\left(\lambda^k, \lambda^* \right) \le C_2 \tau_{EC}^k d\left(\lambda^0, \lambda^* \right)
\end{align}
\end{subequations}
are guaranteed in a contraction region depicted by Riemannian metric
\begin{equation}\label{RiemannianMetric}
M = 
\begin{bmatrix}
\beta c\boldsymbol{I}& \alpha\beta  A_1^\T\\
\alpha\beta  A_1 & \alpha c\boldsymbol{I}
\end{bmatrix}\in\mathbb{R}^{\left(n+p_1\right)\times\left(n+p_1\right)},
\end{equation}
where the convergence rate is 
\begin{equation}\label{ConvergenceRate}
\tau_{EC} = \sqrt{1-\frac{1-c}{c}\alpha\beta\underline{\sigma}_1}.
\end{equation}
with $c = 2\max\left\lbrace\frac{\beta\overline{\sigma}_1}{\underline{\rho}},\alpha\overline{\rho} \right\rbrace$.

\subsection{Inequality Constrained Optimization}

Consider the following inequality constrained optimization problem
\begin{subequations}\label{IECO_OP}
	\begin{align}
	\min_{x\in\mathbb{R}^n}\quad &f\left(x\right)\\
	\text{s.t.} \quad &A_2 x\le b_2 
	\end{align}
\end{subequations}
where $f$ and $A_2$ satisfy Assumptions 1 and 2. Denote $A_2$ and $b_2$ as $A_2 = \left[a_{2,1}^\T,a_{2,2}^\T,\cdots,a_{2,p_2}^\T \right]^\T$ and $b_2=\left[b_{2,1},b_{2,2},\cdots,b_{2,p_2} \right]^\T$, respectively. Note that the algorithm for inequality constrained optimization is based on augmented Lagrangian function. It should be mentioned that there are several constructions of augmented Lagrangian functions in the existing literature \cite{dhingra2018proximal,zhang2018projected,xu2017first}, while we adopt the classic one \cite{xu2017first} in this paper. The augmented Lagrangian of the optimization problem in (\ref{IECO_OP}) is established as follows
\begin{equation}\label{LangrangianIECO}
\mathcal{L}_{IEC}\left(x,\lambda\right) = f\left(x\right)+\varPhi_{\gamma}\left(x, \lambda \right)
\end{equation}
where $\lambda\in\mathbb{R}^{p_2}$ is the Lagrangian multiplier, $\gamma>0$ is a user-defined scalar and $\varPhi_{\gamma}\left(x, \lambda \right) = \sum_{i=1}^{p_2} \varphi_{\gamma}\left( a_{2,i} x-b_{2,i},\lambda_i\right)$ is a penalty function with
\begin{equation}\label{AugLag}
\begin{aligned}
\varphi_{\gamma}&\left( a_{2,i} x-b_{2,i},\lambda_i\right)\\
&= 
\begin{cases}
\left( a_{2,i} x-b_{2,i}\right)\lambda_i+\frac{\gamma}{2}\left(a_{2,i} x-b_{2,i}\right)^2,\\
\qquad\qquad\qquad \text{if } \gamma\left( a_{2,i} x-b_{2,i}\right)+\lambda_i\ge 0;\\
-\frac{\lambda_i^2}{2\gamma},\qquad\quad\;\;\,\text{if } \gamma\left( a_{2,i} x-b_{2,i}\right)+\lambda_i< 0.
\end{cases}
\end{aligned}
\end{equation}
for $i\in\mathbb{N}_{\left[1,p_2\right]}$.

\textit{Remark 2:} Notice that the augmented Lagrangian in (\ref{LangrangianIECO}) used in this paper which is projection free is different from that adopted in \cite{feijer2010stability}. The convergence analysis of the optimization algorithm with projection term in \cite{feijer2010stability} relies on a diagonal Lyapunov function, while the Riemannian metric constructed in this paper is in general form which is less conservative. 

The KKT conditions associated to the problem in (\ref{LangrangianIECO}) for characterizing the optimal pair $\left(x^*,\lambda^* \right)$ are given as
\begin{subequations}\label{KKT2}
	\begin{align}
	\nabla_{x^*} \mathcal{L}_{IEC}\left(x^*,\lambda^*\right) &= \nabla f\left(x^*\right)+ \nabla_x \varPhi_{\gamma}\left(x^*, \lambda^* \right)= 0\\
	\nabla_{\lambda^*} \mathcal{L}_{IEC}\left(x^*,\lambda^*\right) &= \nabla_{\lambda} \varPhi_{\gamma}\left(x^*, \lambda^* \right) = 0
	\end{align}
\end{subequations}
where $\nabla_x \varPhi_{\gamma}\left(x, \lambda \right) = \sum_{i=1}^{p_2} \nabla_x\varphi_{\gamma}\left( a_{2,i}x-b_{2,i},\lambda_i\right)$ and $\nabla_{\lambda} \varPhi_{\gamma}\left(x, \lambda \right) = \sum_{i=1}^{p_2} \nabla_x\varphi_{\gamma}\left( a_{2,i} x-b_{2,i},\lambda_i\right)$ are the gradients with
\begin{subequations}\label{Gradient}
\begin{align}
\nabla_x \varphi_{\gamma}\left(x, \lambda_i \right)& = \left[ \gamma\left( a_{2,i} x-b_{2,i}\right)+\lambda_i \right]_+ a_{2,i}^\T\\
\nabla_{\lambda} \varphi_{\gamma}\left(x, \lambda_i\right) &=\frac{1}{\gamma}\left(\left[  \gamma\left( a_{2,i} x-b_{2,i}\right)+\lambda_i\right]_+-\lambda_i\right) \text{e}_i
\end{align}
\end{subequations}

Thus, the primal-dual gradient dynamics of the optimization updating can be formulated as follows
\begin{subequations}\label{PDGD_IECO}
	\begin{align}
	x^{k+1} &= x^k -\alpha\nabla f\left(x^k\right)\nonumber\\
	&\quad -\alpha \sum_{i=1}^{p_2} \left[ \gamma\left( a_{2,i} x^k-b_{2,i}\right)+\lambda_i^k \right]_+ a_{2,i}^\T\\
	\lambda^{k+1}  &=\lambda^k + \beta \sum_{i=1}^{p_2} \frac{1}{\gamma}\left(\left[  \gamma\left( a_{2,i} x^k-b_{2,i}\right)+\lambda_i^k\right]_+-\lambda_i^k\right) \text{e}_i
	\end{align}
\end{subequations}
where $\alpha$ is the step-size for primal updating and $ \beta $ is for dual updating. To analyze the convergence of the optimization updating dynamics in (\ref{PDGD_IECO}), the following lemma is introduced in addition to Lemma 1. 

\textit{Lemma 3:} Given $\xi = \left[x^\T, \lambda^\T \right]^\T$ and the optimal $\xi^* = \left[\left(x^*\right)^\T, \left(\lambda^*\right)^\T \right]^\T$, there exists a diagonal matrix $\varPsi = \text{diag}\left\lbrace \psi_1, \psi_2, \cdots,\psi_{p_2}\right\rbrace$ where $\psi_i\left(\xi\right)\in\left[0,1\right]$ for $i\in\mathbb{N}_{\left[1,p_2\right]}$ with respect to $\xi$ such that 
\begin{subequations}\label{psi_dynamics}
	\begin{align}
	&\nabla_{x}\phi_{\gamma}\left(x, \lambda_i\right) - \nabla_{x}\phi_{\gamma}\left(x^*, \lambda_i^*\right) \nonumber\\
	&\quad = \gamma\psi_i\left(\xi\right) a_{2,i}\left( x-x^*\right)a_{2,i}^\T+\psi_i\left(\xi\right)\left( \lambda_i-\lambda_i^*\right)a_{2,i}^\T,\\
	&\nabla_{\lambda}\phi_{\gamma}\left(x, \lambda_i\right) - \nabla_{\lambda}\phi_{\gamma}\left(x^*, \lambda_i^*\right) \nonumber\\
	&\quad = \psi_i\left(\xi\right) a_{2,i}\left( x-x^*\right)\text{e}_i+\frac{1}{\gamma}\left(\psi_i\left(\xi\right)-1\right)\left( \lambda_i-\lambda_i^*\right)\text{e}_i.
	\end{align}
\end{subequations}
\textit{Proof:} We define $z_i =\gamma\left( a_{2,i} x-b_{2,i}\right)+\lambda_i$ and $z_i^* =  \gamma\left( a_{2,i} x^*-b_{2,i}\right)+\lambda_i^* $. It can be directly obtained that there exists $\psi_i\in\left[0,1\right]$ satisfying
\begin{equation}\label{psi_dynamics_proof1}
\psi_i = \frac{\left[z_i\right]_+-\left[z_i^*\right]_+}{z_i-z_i^*}.
\end{equation}
The proof is completed.  $\hfill\blacksquare$

The convergence of the primal-dual gradient dynamics in (\ref{PDGD_IECO}), thereafter, is rigorously analyzed and the theoretical analysis results are stated by the following theorem.

\textit{Theorem 2:} Suppose Assumptions 1 and 2 hold. For the primal-dual gradient optimization updating dynamics in (\ref{PDGD_IECO}), if the step-sizes $\alpha$, $\beta$ and the penalty parameter $\gamma$ are properly designed to satisfy $\gamma\ge 2\beta$ and $\max\left\lbrace c_1,c_2,c_3,c_4 \right)\in\left[0,\frac{1}{4}\right] $ where $c_1 = \max\left\lbrace\alpha\overline{\rho},\alpha\gamma\overline{\sigma}_2 \right\rbrace\dot\max\left\lbrace 2, \frac{\overline{\sigma}_2}{\underline{\sigma}_2}\right\rbrace$, $c_2 = 2\beta \frac{\overline{\sigma}_2}{\underline{\sigma}_2}\max\left\lbrace2,\frac{\overline{\sigma}_2}{\underline{\sigma}_2}\right\rbrace$, $c_3 = \frac{\beta\overline{\sigma}_2}{\underline{\rho}}\max\left\lbrace2,\frac{2}{\alpha\gamma\underline{\sigma}_2},\frac{\overline{\sigma}_2}{\alpha\gamma\underline{\sigma}_2^2} \right\rbrace$ and $c_4 = 2\frac{\beta\overline{\sigma}_2^2}{\gamma\underline{\sigma}_2^2}\max\left\lbrace\frac{\overline{\sigma}_2}{\underline{\sigma}_2},\alpha\gamma\overline{\sigma}_2\right\rbrace$, there exist two constants $C_1$ and $C_2$ which depend on $\underline{\rho}$, $\overline{\rho}$, $\underline{\sigma}_1$, $\overline{\sigma}_1$, $\alpha$, $\beta$ and $\gamma$, such that the following conditions
\begin{subequations}\label{IECO_Theorem}
	\begin{align}
	d\left(x^k, x^* \right) \le C_1 \tau_{IEC}^k d\left(x^0, x^* \right)\\
	d\left(\lambda^k, \lambda^* \right) \le C_2 \tau_{IEC}^k d\left(\lambda^0, \lambda^* \right)
	\end{align}
\end{subequations}
are guaranteed in a contraction region depicted by Riemannian metric
\begin{equation}\label{RiemannianMetricIEO}
M = 
\begin{bmatrix}
\beta c\boldsymbol{I}& \alpha\beta  A_2^\T\\
\alpha\beta  A_2 & \alpha c\boldsymbol{I}
\end{bmatrix}\in\mathbb{R}^{\left(n+p_2\right)\times\left(n+p_2\right)},
\end{equation}
where the convergence rate is 
\begin{equation}\label{ConvergenceRateIECO}
\tau_{IEC} = \sqrt{1-\frac{1-c}{c}\alpha\beta\underline{\sigma}_1}.
\end{equation}
with $c = 4\max\left\lbrace c_1,c_2,c_3,c_4 \right\rbrace$.

\section{Numerical Examples}\label{Numerical Examples}

In this section, the presented contraction analysis results are verified by two numerical examples. The equality and inequality constrained optimization cases are simulated, respectively. 

\subsection{Equality Constrained Optimization}

\begin{figure}[htbp]
	\centering{\includegraphics[width=8cm,height=6cm]{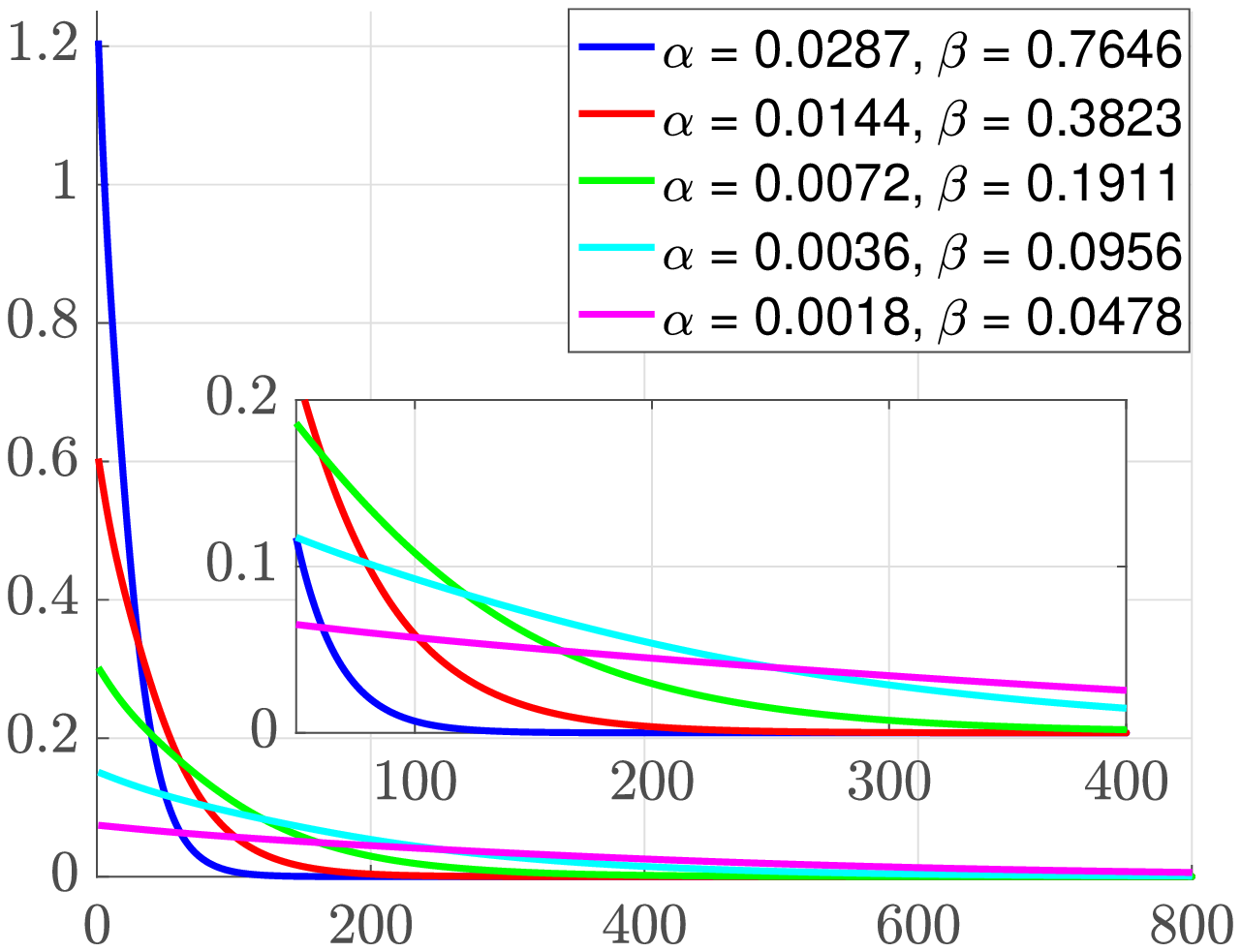}}
	\caption{Riemannian distance between the actual and optimal states.}
	\label{simu:ECO}
\end{figure}

We study a standard quadratic optimization problem. The cost function is $f\left(x\right) = \frac{1}{2}x^\T Q x$ with respect to $x\in\mathbb{R}^{6}$ where $Q = Q_0^\T Q_0 +5\boldsymbol{I}$ with $Q_0$ be a Gaussian random matrix. The equality constraint is depicted by the Gaussian random matrices $A_1\in\mathbb{R}^{3\times 6}$ and $b_1\in\mathbb{R}^{3}$. The step-sizes $\alpha$ and $\beta$ are properly designed according to Theorem 1. Fig. \ref{simu:ECO} shows the Riemannian distance between the actual and optimal states with respect to the Riemannian metric in (\ref{RiemannianMetric}). 

\subsection{Inequality Constrained Optimization}

\begin{figure}[htbp]
	\centering{\includegraphics[width=8cm,height=6cm]{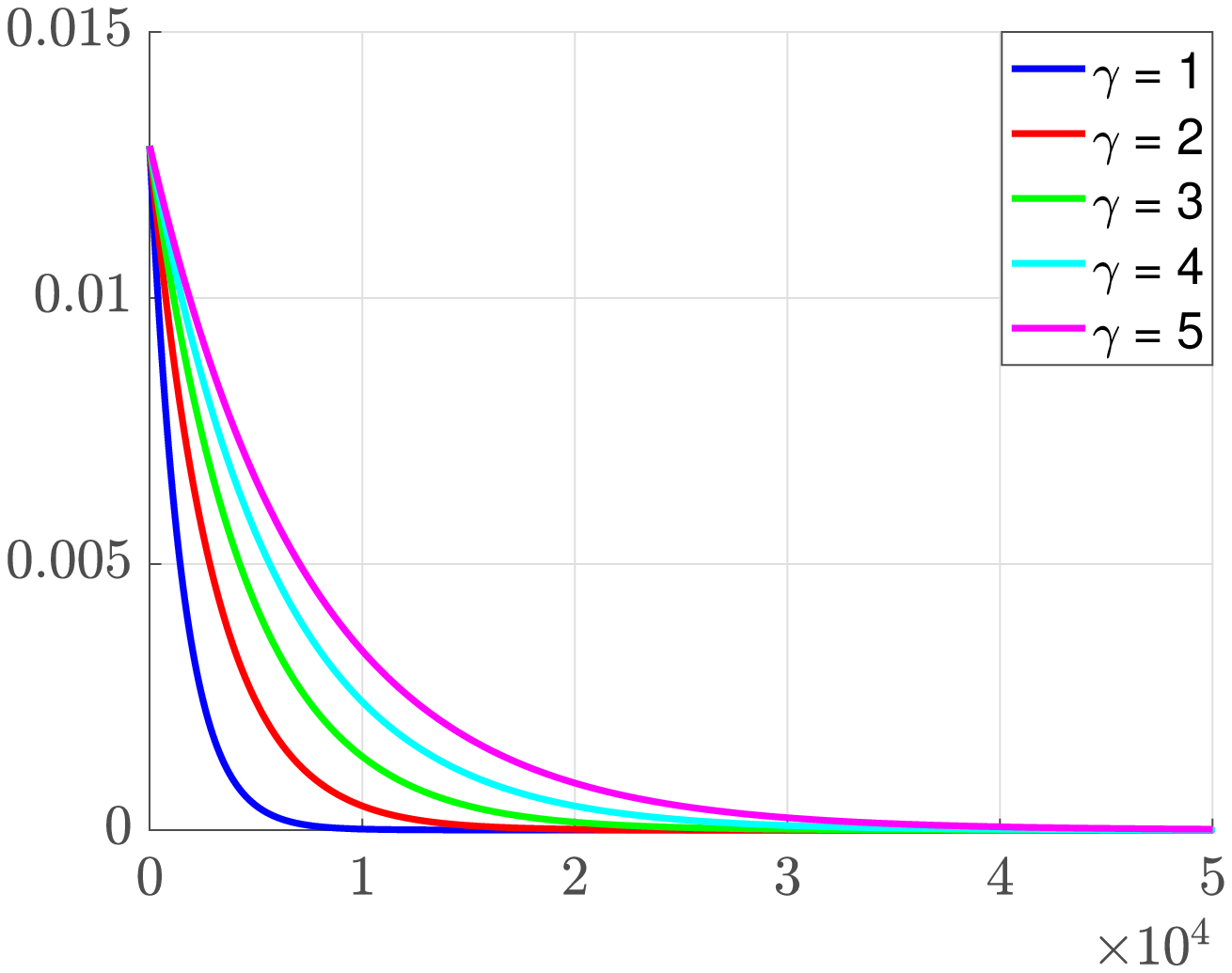}}
	\caption{Riemannian distance between the actual and optimal states.}
	\label{simu:IECO}
\end{figure}

We numerically investigate the similar example as equality constrained optimization. The cost function is $f\left(x\right) = \frac{1}{2}x^\T Q x$ with respect to $x\in\mathbb{R}^{60}$ where $Q = Q_0^\T Q_0 +5\boldsymbol{I}$ with $Q_0$ be a Gaussian random matrix. The inequality constraint is defined by the Gaussian random matrices $A_2\in\mathbb{R}^{30\times 60}$ and $b_2\in\mathbb{R}^{30}$. We carefully design the step-sizes $\alpha$ and $\beta$ according to Theorem 2. Fig. \ref{simu:IECO} compares the Riemannian distances between the actual and optimal states with respect to the Riemannian metric in (\ref{RiemannianMetric}) under different penalty parameter $\gamma$.

\section{Conclusion}\label{Conclusion}

In this paper, we investigate the primal-dual gradient optimization in the context of discrete-time updating dynamics. The contraction theory based on Riemannian manifolds is first utilized to analyze the geometric convergence of a convex optimization algorithm. The equality and inequality constrained optimization cases are studied, respectively. Under some rational assumptions, the Riemannian metric is constructed to span a contraction region. It is shown that if the step-sizes of the updating dynamics are properly designed, the convergence rates for both cases can be obtained based on the contraction region in which the convergence can be guaranteed. Furthermore, we adopt the augmented Lagrangian function which is projection free to handle the inequality constraints. Two numerical experiments are simulated to demonstrate the effectiveness of the proposed primal-dual gradient optimization algorithm.


%

\appendices
\section{Proof of Theorem 1}

Stack $x$ and $\lambda$ into a vector $\xi = \left[x^\T, \lambda^\T \right]^\T$. For the optimal pair, it can be similarly defined as $\xi^* = \left[ \left(x^* \right)^\T,\left(\lambda^* \right)^\T \right]$. In the following proof, we drop the dependence of $G\left(x\right)$ for simplification. According to the primal-dual gradient updating dynamics in (\ref{PDGD_ECO}), we can obtain the differential dynamics as follows
\begin{equation}\label{ECO_Theorem_Proof1}
\begin{aligned}
&\xi^+-\left(\xi^*\right)^+\\
&= 
\begin{bmatrix}
x^k-\alpha \nabla f\left(x^k\right)-\alpha A_1^\T \lambda^k -x^* +\alpha\nabla f\left(x^*\right)-\alpha A_1^\T\lambda^*\\
\lambda^k+\beta \left(A_1 x^k-b_1 \right)-\lambda^*+\beta \left(A_1 x^*-b_1 \right)
\end{bmatrix}\\
&=
\begin{bmatrix}
\left(x^k-x^*\right)-\alpha G\left( x^k-x^*\right)- \alpha A_1^\T \left(\lambda^k -\lambda^*\right)\\
\beta A_1\left( x^k-x^*\right)+\left(\lambda^k-\lambda^*\right)
\end{bmatrix}\\
&\triangleq\Theta\left(\xi - \xi^* \right)
\end{aligned}
\end{equation}
where 
\begin{equation}\label{Theta}
\Theta = \begin{bmatrix}
\boldsymbol{I}-\alpha G & -\alpha A_1^\T \\
\beta A_1 & \boldsymbol{I}
\end{bmatrix}.
\end{equation}
By constructing the Riemannian metric in (\ref{RiemannianMetric}), the difference of the Riemannian energy between the adjacent updates can be written as
\begin{equation}\label{diffRieDist}
\begin{aligned}
E\left(\xi^+,\left(\xi^*\right)^+\right)&-E\left(\xi,\xi^*\right)\\
=&\left\| \xi^+-\left(\xi^*\right)^+\right\|_{M}^2-\left\| \xi-\xi^*\right\|_{M}^2\\
=&\left(\xi-\xi^*\right)^\T \left(\Theta^\T M\Theta - M  \right) \left(\xi-\xi^*\right)\\
\triangleq&\left(\xi-\xi^*\right)^\T\Pi\left(\xi-\xi^*\right),
\end{aligned}
\end{equation}
where 
\begin{equation}\label{Pi}\nonumber
\Pi = 
\begin{bmatrix}
\Pi_{1} & \Pi_{2}^\T\\
\Pi_{2} & \Pi_{3}
\end{bmatrix}
\end{equation}
with
\begin{subequations}
	\begin{align}
	\Pi_{1} &= \beta c\left(\boldsymbol{I}-\alpha G \right)^\T\left(\boldsymbol{I}-\alpha G \right)+\alpha\beta^2 \left(A_1^\T A_1\right.\nonumber\\
	&\quad \left.\times\left(\boldsymbol{I}-\alpha G\right)+\left(\boldsymbol{I}-\alpha G\right)^\T A_1^\T A_1\right)\nonumber\\
	&\quad +\alpha\beta^2 c A_1^\T A_1 -\beta c \boldsymbol{I},\\
	\Pi_{2} &=  \alpha\beta\left(1-c\right)A_1\left(\boldsymbol{I}-\alpha G \right)-\alpha^2 \beta^2 A_1 A_1^\T A_1\nonumber\\
	&\quad - \alpha\beta\left(1-c\right) A_1\\
	\Pi_{3} &= \left(\alpha^2\beta c-2\alpha^2\beta \right) A_1 A_1^\T .
	\end{align}
\end{subequations}
In what follows, we need to guarantee the contraction property in the contraction region depicted by $M$. To prove $\left(\ref{contractionTheoryDef} \right)$, it suffices to prove that $\Pi\le\left(\tau_{EC}^2-1\right) M$. Letting $\varPhi = \left(\tau_{EC}^2-1\right) M-\Pi$, we can obtain
\begin{equation}\label{varPhi}\nonumber
\varPhi = 
\begin{bmatrix}
\varPhi_{1}&\varPhi_{2}^\T\\
\varPhi_{2}&\varPhi_{3}
\end{bmatrix}
\end{equation}
with
\begin{subequations}
	\begin{align}
	\varPhi_{1} &= \tau_{EC}^2\beta c\boldsymbol{I}- \beta c\left(\boldsymbol{I}-\alpha G \right)^\T\left(\boldsymbol{I}-\alpha G \right)-\alpha\beta^2 \left(A_1^\T A_1\right.\nonumber\\
	&\quad \times\left.\left(\boldsymbol{I}-\alpha G\right)+\left(\boldsymbol{I}-\alpha G\right)^\T A_1^\T A_1\right)\nonumber\\
	&\quad -\alpha\beta^2 c A_1^\T A_1,\\
	\varPhi_{2} &=  \tau_{EC}^2\alpha\beta A_1 -\alpha\beta\left(1-c\right)A_1\left(\boldsymbol{I}-\alpha G \right)\nonumber\\
	&\quad +\alpha^2 \beta^2 A_1 A_1^\T A_1 - \alpha\beta c A_1\\
	\varPhi_{3} &= \left(\tau_{EC}^2-1\right)\alpha c\boldsymbol{I}-\left(\alpha^2\beta c-2\alpha^2\beta \right) A_1 A_1^\T.
	\end{align}
\end{subequations}
Thus, it is to show $\varPhi\ge 0$ in the following. Resorting to the Schur's complement, to prove $\varPhi\ge0$, it is sufficient to prove $\varPhi_{3}>0$ and $\varPhi_{1}-\varPhi_{2}\varPhi_{3}^{-1}\varPhi_{2}^\T\ge 0$. 

Under Assumption 1, by recalling the convergence rate in (\ref{ConvergenceRate}), we can get
\begin{equation}\label{varPhi3}
\begin{aligned}
\varPhi_{3} &= \left(1-\frac{1-c}{c}\alpha \beta \underline{\sigma}_1-1\right)\alpha c\boldsymbol{I}-\left(\alpha^2\beta c-2\alpha^2\beta \right) A_1 A_1^\T \\
&\ge \alpha^2 \beta A_1 A_1^\T> 0
\end{aligned}
\end{equation}
Moreover, in accordance with $A^\T \left(A A^\T \right)^{-1} A\le \boldsymbol{I}$, we can get
\begin{equation}\label{ECO_Theorem_Proof2}
\begin{aligned}
&\varPhi_{2}^\T\varPhi_{3}^{-1}\varPhi_{2}\\
& \le \varPhi_{2}\left(\alpha^2 \beta A_1 A_1^\T\right)^{-1}\varPhi_{2}^\T\\
& = \beta \left(\left(\tau_{EC}^2-1\right)\boldsymbol{I}+\alpha\left(1-c\right)G\right)^\T\left(\left(\tau_{EC}^2-1\right)\boldsymbol{I}\right.\\
& \quad \left.+\alpha\left(1-c\right)G\right)+\alpha^2\beta^3 A_1^\T A_1 A_1^\T A_1\\
& \quad +\alpha\beta^2\left(\left(\left(\tau_{EC}^2-1\right)\boldsymbol{I}+\alpha\left(1-c\right)G\right)^\T A_1^\T A_1\right.\\
& \quad \left. + A_1^\T A_1\left(\left(\tau_{EC}^2-1\right)\boldsymbol{I}+\alpha\left(1-c\right)G\right)\right)\\
& \le \frac{1-4c+3c^2}{c^2}\alpha^2\beta^3\underline{\sigma}_1^2\boldsymbol{I}+\alpha^2\beta\left(1-c\right)^2 \overline{\rho} G\\
&\quad +\alpha^2\beta^3\overline{\sigma}_1^2\boldsymbol{I}-\frac{\left(2c-1\right)^2+1}{c}\alpha^2\beta^2\underline{\sigma}_1 G\\
&\quad +\alpha^2\beta^2\left(G^\T A_1^\T A_1+A_1^\T A_1 G\right)\\
\end{aligned}
\end{equation}
Thus, recalling that $c = 2\max\left\lbrace\frac{\beta\overline{\sigma}_1}{\underline{\rho}},\alpha\overline{\rho} \right\rbrace\in \left[0,1\right] $, one can have
\begin{equation}\label{ECO_Theorem_Proof3}
\begin{aligned}
&\varPhi_{1}-\varPhi_{2}^\T\varPhi_{3}^{-1}\varPhi_{2}\\
& \ge -\alpha\beta^2\left(1-c\right)\underline{\sigma}_1\boldsymbol{I}+2\alpha\beta cG -\alpha^2\beta c \overline{\rho} G\\
& \quad -2\alpha\beta^2 \overline{\sigma}_1\boldsymbol{I} -\alpha\beta^2 c \overline{\sigma}_1\boldsymbol{I}-\alpha^2\beta\left(1-c\right)^2 \overline{\rho} G\\
&\quad -\frac{1-4c+3c^2}{c^2}\alpha^2\beta^3\underline{\sigma}_1^2\boldsymbol{I}-\alpha^2\beta^3\overline{\sigma}_1^2\boldsymbol{I}\\
&\quad +\frac{\left(2c-1\right)^2+1}{c}\alpha^2\beta^2\underline{\sigma}_1 G\\
&\ge 2\alpha\beta c G+\frac{\left(2c-1\right)^2+1}{c}\alpha^2\beta^2\underline{\sigma}_1\overline{\rho}\boldsymbol{I}-\alpha\beta^2\underline{\sigma}_1\boldsymbol{I}\\
&\quad -2\alpha\beta^2\overline{\sigma}_1\boldsymbol{I}-\frac{\left(2c-1\right)^2}{c^2}\alpha^2\beta^3\overline{\sigma}_1^2\boldsymbol{I}\\
&\quad -\alpha^2\beta\left(\frac{3}{4}+\left(c-\frac{1}{2}\right)^2\right)\overline{\rho} G\ge 0
\end{aligned}
\end{equation}
where the last inequality follows from: i) $\frac{3}{2}\alpha\beta c G\ge 3\alpha\beta^2\overline{\sigma}_1 \boldsymbol{I}\ge \alpha\beta^2\underline{\sigma}_1\boldsymbol{I}+2\alpha\beta^2\overline{\sigma}_1\boldsymbol{I}$, ii) $\frac{1}{2}\alpha\beta c G \ge \alpha^2\beta\overline{\rho} G\ge\alpha^2\beta\left(\frac{3}{4}+\left(c-\frac{1}{2}\right)^2\right)\overline{\rho} G$, and iii) $\frac{\left(2c-1\right)^2+1}{c}\alpha^2\beta^2\underline{\sigma}_1\overline{\rho}\boldsymbol{I}\ge \frac{\left(2c-1\right)^2}{c^2}\alpha^2\beta^3\overline{\sigma}_1^2\boldsymbol{I}$.

The geometric convergence can be guaranteed within the contraction region, which implies Eq. (\ref{ECO_Theorem}) is satisfied. The proof is completed. $\hfill\blacksquare$

\section{Proof of Theorem 2}
To prove Theorem 2, we stack $x$ and $\lambda$ into a vector $\xi = \left[x^\T, \lambda^\T \right]^\T$. For the optimal pair, we adopt the similarly definition $\xi^* = \left[ \left(x^* \right)^\T,\left(\lambda^* \right)^\T \right]$. Inspired by Lemma 3, we can obtain the differential dynamics as follows according to the primal-dual gradient updating dynamics in (\ref{PDGD_IECO})
\begin{equation}\label{IECO_Theorem_Proof1}
\begin{aligned}
x^+&-\left(x^*\right)^+\\
=& x^k-\alpha \nabla f\left(x^k\right)-\alpha \sum_{i=1}^{p_2} \left[ \gamma\left( a_{2,i} x^k-b_{2,i}\right)+\lambda_i^k \right]_+ a_{2,i}^\T \\
&-x^* +\alpha\nabla f\left(x^*\right)-\alpha \sum_{i=1}^{p_2} \left[ \gamma\left( a_{2,i} x^*-b_{2,i}\right)+\lambda_i^* \right]_+ a_{2,i}^\T\\
=&\left(x^k-x^*\right)-\alpha\left(\nabla f\left(x^k\right)-f\left(x^*\right)\right)-\alpha\sum_{i=1}^{p_2}\left(\left[ \gamma\left( a_{2,i} x^k\right.\right.\right.\\
&\left.\left.\left.-b_{2,i}\right)+\lambda_i^k \right]_+ -\left[ \gamma\left( a_{2,i} x^*-b_{2,i}\right)+\lambda_i^* \right]_+  \right)a_{2,i}^\T\\
=&\left(\boldsymbol{I}-\alpha G-\alpha A_2^\T\varPsi A_2\right)\left(x^k-x^*\right)-\alpha A_2^\T\varPsi\left(\lambda^k-\lambda^*\right)
\end{aligned}
\end{equation}
and
\begin{equation}\label{IECO_Theorem_Proof2}
\begin{aligned}
\lambda^+&-\left(\lambda^*\right)^+\\
=& \lambda^k+ \beta \sum_{i=1}^{p_2} \frac{1}{\gamma}\left(\left[  \gamma\left( a_{2,i} x^k-b_{2,i}\right)+\lambda_i^k\right]_+-\lambda_i^k\right) \text{e}_i \\
&-\lambda^* - \beta \sum_{i=1}^{p_2} \frac{1}{\gamma}\left(\left[  \gamma\left( a_{2,i} x^*-b_{2,i}\right)+\lambda_i^*\right]_+-\lambda_i^*\right) \text{e}_i\\
=&\beta \varPsi A_2\left(x^k-x^*\right)+\left(\boldsymbol{I}-\frac{\beta}{\gamma}\left(\boldsymbol{I}-\varPsi \right) \right)\left(\lambda^k-\lambda^*\right)
\end{aligned}
\end{equation}
Note that we drop the dependence of $\varPsi\left(\xi\right)$ in (\ref{IECO_Theorem_Proof1}) and (\ref{IECO_Theorem_Proof2}) for simplification. Thus, one has
\begin{equation}\label{IECO_Theorem_Proof3}
\xi^+-\left(\xi^*\right)^+\triangleq\Theta\left(\xi - \xi^* \right)
\end{equation}
where 
\begin{equation}\label{Theta2}
\Theta = \begin{bmatrix}
\boldsymbol{I}-\alpha G-\alpha A_2^\T\varPsi A_2 & -\alpha A_2^\T\varPsi \\
\beta \varPsi A_2 & \boldsymbol{I}-\frac{\beta}{\gamma}\left(\boldsymbol{I}-\varPsi\right)
\end{bmatrix}.
\end{equation}
By constructing the Riemannian metric in (\ref{RiemannianMetricIEO}), consider the difference of the Riemannian energy between the adjacent updates as follows
\begin{equation}\label{diffRieDist_IECO}
\begin{aligned}
E\left(\xi^+,\left(\xi^*\right)^+\right)&-E\left(\xi,\xi^*\right)\\
=&\left\| \xi^+-\left(\xi^*\right)^+\right\|_{M}^2-\left\| \xi-\xi^*\right\|_{M}^2\\
=&\left(\xi-\xi^*\right)^\T \left(\Theta^\T M\Theta - M  \right) \left(\xi-\xi^*\right)\\
\triangleq&\left(\xi-\xi^*\right)^\T\Pi\left(\xi-\xi^*\right),
\end{aligned}
\end{equation}
where 
\begin{equation}\label{Pi_IECO}\nonumber
\Pi = 
\begin{bmatrix}
\Pi_{1} & \Pi_{2}^\T\\
\Pi_{2} & \Pi_{3}
\end{bmatrix}
\end{equation}
with
\begin{subequations}
	\begin{align}
	\Pi_{1} &= \beta c\left(\boldsymbol{I}-\alpha G-\alpha A_2^\T\varPsi A_2 \right)^\T\left(\boldsymbol{I}-\alpha G -\alpha A_2^\T\varPsi A_2\right)\nonumber\\
	&\quad +\alpha\beta^2 \left(A_2^\T \varPsi A_2\left(\boldsymbol{I}-\alpha G-\alpha A_2^\T\varPsi A_2\right)+\left(\boldsymbol{I}-\alpha G\right.\right.\nonumber\\
	&\left.\left.\quad -\alpha A_2^\T\varPsi A_2\right)^\T A_2^\T \varPsi A_2\right) +\alpha\beta^2 c A_2^\T \varPsi \varPsi A_2 -\beta c \boldsymbol{I},\\
	\Pi_{2} &= -A_2\left(\alpha^2\beta \left(1-c\right)\left(G+ \gamma A_2^\T\varPsi A_2\right)+\alpha^2\beta^2 A_2^\T\varPsi A_2\right)\nonumber \\
	&\quad -\left(\boldsymbol{I}-\varPsi\right)\left(\alpha^2\beta \left(c-\frac{\beta}{\gamma}\right)A_2\left(G+ \gamma A_2^\T\varPsi A_2\right)\right.\nonumber\\
	&\quad \left. -\alpha^2\beta^2 A_2 A_2^\T\varPsi A_2 -\frac{\alpha\beta^2 c}{\gamma} \left(\boldsymbol{I}-\varPsi\right)A_2 \right.\nonumber\\
	&\quad \left.+\frac{\alpha\beta^2 \left(c+1\right)}{\gamma}A_2\right)\\
	\Pi_{3} &= \alpha^2\beta\left(\left(c-\frac{2\beta}{\gamma}\right)\varPsi A_2 A_2^\T\varPsi -\left(1-\frac{\beta}{\gamma} \right)\left( A_2 A_2^\T\varPsi\right.\right.\nonumber\\
	&\quad \left.\left.+\varPsi A_2 A_2^\T \right) \right)+\alpha c \left(\frac{\beta^2}{\gamma^2}\left(\boldsymbol{I}-\varPsi \right)-\frac{2\beta}{\gamma}\boldsymbol{I} \right)\left(\boldsymbol{I}-\varPsi \right). \label{Pi_3_IECO}
	\end{align}
\end{subequations}
The contraction property can be guaranteed within the contraction region. Therefore, to prove Eq. $\left(\ref{contractionTheoryDef} \right)$, it suffices to prove that $\Pi\le\left(\tau_{IEC}^2-1\right) M$. Letting $\varPhi = \left(\tau_{IEC}^2-1\right) M-\Pi$, we can obtain
\begin{equation}\label{varPhi_IECO}\nonumber
\varPhi = 
\begin{bmatrix}
\varPhi_{1}&\varPhi_{2}^\T\\
\varPhi_{2}&\varPhi_{3}
\end{bmatrix}
\end{equation}
with
\begin{subequations}
	\begin{align}
	\varPhi_{1} &= \left(\tau_{IEC}^2-1\right)\beta c\boldsymbol{I}-\Pi_{1},\\
	\varPhi_{2} &=  \left(\tau_{IEC}^2-1\right)\alpha\beta A_2 -\Pi_{2},\\
	\varPhi_{3} &= \left(\tau_{IEC}^2-1\right)\alpha c\boldsymbol{I}-\Pi_{3}.
	\end{align}
\end{subequations}
Thus, to prove Eq. $\left(\ref{contractionTheoryDef} \right)$, it is sufficient to show $\varPhi\ge 0$. Resorting to the Schur's complement, to prove $\varPhi\ge0$, it suffices to prove $\varPhi_{3}>0$ and $\varPhi_{1}-\varPhi_{2}\varPhi_{3}^{-1}\varPhi_{2}^\T\ge 0$, respectively.

We study the upper bound of $\Pi_{3}$ to bound $\varPhi_{3}$ and summarize the analysis results in Lemma 4. The proof of Lemma 4 is deferred to Appendix C.

\textit{Lemma 4:} For the inequality constrained optimization problem in (\ref{IECO_OP}), if it holds that $\varPsi = \text{diag}\left(\psi_1,\psi_2,\cdots,\psi_{p_2} \right)$ with $\psi_i\in\left[0,1\right]$ for $i\in\mathbb{N}_{\left[1,m\right]}$, $\gamma\ge 2\beta$ and
\begin{equation}\label{Pi_3_IECO_bound1}
\frac{2\alpha\beta\left(1-\frac{\beta}{\gamma}\right)^2\overline{\sigma}_2}{1-\left(1-\frac{\beta}{\gamma}\right)^2}\le c \le 1,
\end{equation}
the matrix $\Pi_{3}$ defined in (\ref{Pi_3_IECO}) satisfies
\begin{equation}\label{Pi_3_IECO_bound2}
\Pi_{3} \le -\frac{1}{2}\alpha^2\beta c A_2 A_2^\T.
\end{equation}

Inspired by Lemma 4, under Assumptions 1 and 2, we can obtain the following expression by recalling the convergence rate in (\ref{ConvergenceRateIECO})
\begin{equation}\label{varPhi3_IECO}
\begin{aligned}
\varPhi_{3} &= \left(1-\frac{1-c}{c}\alpha \beta \underline{\sigma}_2-1\right)\alpha c\boldsymbol{I}-\Pi_{3}\\
&\ge -\left(1-c\right)\alpha^2 \beta \underline{\sigma}_2+\left(\frac{3}{2}-c\right)\alpha^2\beta A_2 A_2^\T\\
&\ge \frac{1}{2}\alpha^2\beta A_2 A_2^\T >0
\end{aligned}
\end{equation}
where the first inequality follows from Lemma 4.

Furthermore, rewrite $\varPhi_{2}$ into $\varPhi_{2} = \alpha\beta A_2\varPhi_2^1 + \alpha\beta \left(\boldsymbol{I}-\varPsi\right)\varPhi_2^2$, where
\begin{subequations}\label{Pi_2_rewrite}
	\begin{align}
	\varPhi_2^1 &= \alpha \left(1-c\right)\left(G+ \gamma A_2^\T\varPsi A_2\right)+\alpha\beta A_2^\T\varPsi A_2-\frac{1-c}{c}\alpha\beta\underline{\sigma}_2\boldsymbol{I}\\
	\varPhi_2^2 &= \alpha \left(c-\frac{\beta}{\gamma}\right)A_2\left(G+ \gamma A_2^\T\varPsi A_2\right) -\alpha\beta A_2 A_2^\T\varPsi A_2 \nonumber\\
	& \quad+\frac{\beta}{\gamma} \left(\boldsymbol{I}+c\varPsi\right) A_2.
	\end{align}
\end{subequations}
Thus, one can have
\begin{equation}\label{upperBounding}
\begin{aligned}
\varPhi_{2}^\T \varPhi_{3}^{-1} \varPhi_{2}&\le \frac{2}{\alpha^2\beta}\varPhi_{2}^\T\left(A_2 A_2^\T\right)^{-1}  \varPhi_{2}\\
&\le2\beta\left(\left(\varPhi_2^1\right)^\T \varPhi_2^1 + \left\|\left(\varPhi_2^2\right)^\T\left(A_2 A_2^\T\right)^{-1} \varPhi_2^2\right\|\boldsymbol{I}\right.\\
&\quad \left.+ 2\left\|\left(\varPhi_2^1\right)^\T A_2^\T\left(A_2 A_2^\T\right)^{-1} \varPhi_2^2\right\|\boldsymbol{I}\right)
\end{aligned}
\end{equation}
To further bound $\varPhi_{2}^\T \varPhi_{3}^{-1} \varPhi_{2}$, we utilize Eq. (\ref{upperBounding_1})-Eq. (\ref{upperBounding_3}),
\begin{figure*}
	\begin{equation}\label{upperBounding_1}
	\begin{aligned}
	&\left(\varPhi_2^1\right)^\T \varPhi_2^1 \le \frac{\overline{\sigma}_2}{\underline{\sigma}_2}\left\|\varPhi_2^1\right\|^2 \le \frac{\overline{\sigma}_2}{\underline{\sigma}_2}\left( \alpha\left(1-c\right)\left(\overline{\rho}+\gamma \overline{\sigma}_2\right)+\left(\alpha\beta\overline{\sigma}_2+\frac{1-c}{c}\alpha\beta\underline{\sigma}_2\right) \right)^2\\
	& = \frac{\overline{\sigma}_2}{\underline{\sigma}_2}\left(\alpha^2\left(1-c\right)^2\left(\overline{\rho}+\gamma \overline{\sigma}_2\right)^2+2\alpha\left(1-c\right)\left(\alpha\beta\overline{\sigma}_2+\frac{1-c}{c}\alpha\beta\underline{\sigma}_2\right)\left(\overline{\rho}+\gamma \overline{\sigma}_2\right)+\left(\alpha\beta\overline{\sigma}_2+\frac{1-c}{c}\alpha\beta\underline{\sigma}_2\right)^2\right)
	\end{aligned}
	\end{equation}
	\hrule
	\vspace{3pt}
	\begin{equation}\label{upperBounding_2}
	\begin{aligned}
	&\left\|\left(\varPhi_2^2\right)^\T\left(A_2 A_2^\T\right)^{-1} \varPhi_2^2\right\| \le  \left\| \varPhi_2^2 \right\|^2 \left\| \left(A_2 A_2^\T\right)^{-1} \right\|\le \frac{\overline{\sigma}_2}{\underline{\sigma}_2} \left(\alpha \left(c-\frac{\beta}{\gamma}\right) \left(\overline{\rho}+\gamma\overline{\sigma}_2 \right)+\alpha\beta\overline{\sigma}_2+\frac{\beta \left(1+c\right)}{\gamma}\right)^2\\
	& = \frac{\overline{\sigma}_2}{\underline{\sigma}_2} \left(\alpha^2 \left(c-\frac{\beta}{\gamma}\right)^2 \left(\overline{\rho}+\gamma\overline{\sigma}_2 \right)^2+2\alpha \left(c-\frac{\beta}{\gamma}\right)\left(\frac{\beta \left(1+c\right)}{\gamma}+\alpha\beta\overline{\sigma}_2\right)\left(\overline{\rho}+\gamma\overline{\sigma}_2 \right)+\left(\frac{\beta \left(1+c\right)}{\gamma}+\alpha\beta\overline{\sigma}_2\right)^2\right)
	\end{aligned}
	\end{equation}
	\hrule
	\vspace{3pt}
	\begin{equation}\label{upperBounding_3}
	\begin{aligned}
	&\left\|\left(\varPhi_2^1\right)^\T A_2^\T\left(A_2 A_2^\T\right)^{-1} \varPhi_2^2\right\|\le  \left\| \varPhi_2^1 \right\| \left\| A_2^\T \right\| \left\| \left(A_2 A_2^\T\right)^{-1} \right\| \left\| \varPhi_2^2 \right\|\\
	& \le \frac{\overline{\sigma}_2}{\underline{\sigma}_2}  \left( \alpha\left(1-c\right)\left(\overline{\rho}+\gamma \overline{\sigma}_2\right)+\alpha\beta\overline{\sigma}_2 +\frac{1-c}{c}\alpha\beta\underline{\sigma}_2\right)\left(\alpha \left(c-\frac{\beta}{\gamma}\right) \left(\overline{\rho}+\gamma\overline{\sigma}_2 \right)+\alpha\beta\overline{\sigma}_2+\frac{\beta \left(1+c\right)}{\gamma}\right)\\
	& = \frac{\overline{\sigma}_2}{\underline{\sigma}_2}  \left( \alpha^2\left(1-c\right)\left(c-\frac{\beta}{\gamma}\right)\left(\overline{\rho}+\gamma \overline{\sigma}_2\right)^2+\alpha\left(\left(c-\frac{\beta}{\gamma} \right)\left(\alpha\beta\overline{\sigma}_2 +\frac{1-c}{c}\alpha\beta\underline{\sigma}_2\right)\right.\right.\\
	&\quad\left.\left.+\left(1-c\right)\left(\frac{\beta \left(1+c\right)}{\gamma}+\alpha\beta\overline{\sigma}_2\right)\right)\left(\overline{\rho}+\gamma \overline{\sigma}_2\right)+\left(\alpha\beta\overline{\sigma}_2 +\frac{1-c}{c}\alpha\beta\underline{\sigma}_2\right)\left(\frac{\beta \left(1+c\right)}{\gamma}+\alpha\beta\overline{\sigma}_2\right) \right)
	\end{aligned}
	\end{equation}
	\hrule
	\vspace{3pt}
	\begin{equation}\label{lowerBounding}
	\begin{aligned}
	\varPhi_{1}-\varPhi_{2}^\T \varPhi_{3}^{-1} \varPhi_{2}&\ge -\left(1-c\right)\alpha\beta^2\underline{\sigma}_2 - \alpha\beta^2\left(2+c\right)\overline{\sigma}_2+\left(2\alpha\beta c+2\alpha^2\beta^2\underline{\sigma}_2 \right)\left(G+\gamma A_2^\T \varPsi A_2 \right)\\
	&\quad-\alpha^2\beta c\left(\overline{\rho}+\gamma\overline{\sigma}_2\right)\left(G+\gamma A_2^\T \varPsi A_2\right)-\varPhi_{2}^\T \varPhi_{3}^{-1} \varPhi_{2}\\
	&\ge - 3\alpha\beta^2\overline{\sigma}_2+\left(2\alpha\beta c+2\alpha^2\beta^2\underline{\sigma}_2 \right)\left(\underline{\rho}+\gamma \underline{\sigma}_2 \right) -\alpha^2\beta c\left(\overline{\rho}+\gamma\overline{\sigma}_2\right)\left(G+\gamma A_2^\T \varPsi A_2\right)\\
	&\quad -\frac{2\beta\overline{\sigma}_2}{\underline{\sigma}_2}\left(\alpha^2\left(1-\frac{\beta}{\gamma}\right)^2\left(\overline{\rho}+\gamma\overline{\sigma}_2 \right)^2 +2\alpha\beta\left(1-\frac{\beta}{\gamma}\right)\left(\frac{\beta}{\gamma}+\frac{\beta}{\gamma}c+\alpha\beta\overline{\sigma}_2+\frac{\alpha\beta\overline{\sigma}_2}{c}\right)\left(\overline{\rho}+\gamma\overline{\sigma}_2 \right)\right.\\
	&\quad\left.\frac{\beta^2}{\gamma^2}+\frac{2\beta^2 c}{\gamma^2}+\frac{\beta^2 c^2}{\gamma^2}+\frac{2\alpha\beta^2\overline{\sigma}_2}{\gamma c} +\frac{4\alpha\beta^2\overline{\sigma}_2}{\gamma}+\frac{2\alpha\beta^2\overline{\sigma}_2 c}{\gamma}+\frac{\alpha^2\beta^2\overline{\sigma}_2^2}{c^2} +\frac{2\alpha^2\beta^2\overline{\sigma}_2^2}{c}+\alpha^2\beta^2\overline{\sigma}_2^2\right)\ge 0
	\end{aligned}
	\end{equation}
	\hrule
\end{figure*}
where the last inequality in Eq. (\ref{upperBounding_3}) follows from: i) $\frac{3}{8}\alpha\beta c\underline{\rho}\ge 3\alpha\beta^2\overline{\sigma}$; ii) $\frac{1}{8}\alpha\beta c\ge\frac{1}{8}\alpha\beta c^2\ge \alpha^2\beta c\left(\overline{\rho}+\gamma\overline{\sigma}_2\right)$; iii) $\frac{1}{2}\left(1-\frac{\beta}{\gamma}\right)^2\alpha\beta c\ge\frac{2\beta\overline{\sigma}_2}{\underline{\sigma}_2}\alpha^2\left(1-\frac{\beta}{\gamma}\right)^2\left(\overline{\rho}+\gamma\overline{\sigma}_2\right)$; iv) $\frac{1}{4}\frac{\beta}{\gamma}\left(1-\frac{\beta}{\gamma}\right)\alpha\beta c\ge2\alpha\beta \frac{2\beta\overline{\sigma}_2}{\underline{\sigma}_2}\left(1-\frac{\beta}{\gamma}\right)\frac{\beta}{\gamma}$; v) $\frac{1}{4}\frac{\beta}{\gamma}\left(1-\frac{\beta}{\gamma}\right)\alpha\beta c\ge\frac{1}{4}\frac{\beta}{\gamma}\left(1-\frac{\beta}{\gamma}\right)\alpha\beta c^2\ge2\alpha\beta \frac{2\beta\overline{\sigma}_2}{\underline{\sigma}_2}\left(1-\frac{\beta}{\gamma}\right)\frac{\beta}{\gamma}c$; vi) $\frac{1}{2}\left(1-\frac{\beta}{\gamma}\right)\alpha^2\beta^2\underline{\sigma}_2 \ge2\alpha\beta\frac{2\beta\overline{\sigma}_2}{\underline{\sigma}_2}\left(1-\frac{\beta}{\gamma}\right)\frac{2\alpha\beta^2\overline{\sigma}_2}{\gamma c}$; vii) $\frac{1}{2}\left(1-\frac{\beta}{\gamma}\right)\alpha\beta c\ge2\alpha\beta\frac{2\beta\overline{\sigma}_2}{\underline{\sigma}_2}\left(1-\frac{\beta}{\gamma}\right)2\alpha\beta^2\overline{\sigma}_2$; viii) $\frac{1}{4}\frac{\beta}{\gamma}\alpha\beta c \underline{\rho}\ge\frac{2\beta\overline{\sigma}_2}{\underline{\sigma}_2}\frac{\beta^2}{\gamma^2} $; ix) $\frac{1}{2}\frac{\beta}{\gamma}\alpha\beta c \underline{\rho}\ge\frac{1}{2}\frac{\beta}{\gamma}\alpha\beta c^2 \underline{\rho}\ge\frac{2\beta\overline{\sigma}_2}{\underline{\sigma}_2}\frac{2\beta^2c}{\gamma^2} $; x) $\frac{1}{4}\frac{\beta}{\gamma}\alpha\beta c \underline{\rho}\ge\frac{1}{4}\frac{\beta}{\gamma}\alpha\beta c^3 \underline{\rho}\ge\frac{2\beta\overline{\sigma}_2}{\underline{\sigma}_2}\frac{\beta^2c^2}{\gamma^2}$; xi) $\frac{1}{2}\alpha^2\beta^2\underline{\sigma}_2\underline{\rho}\ge\frac{2\beta\overline{\sigma}_2}{\underline{\sigma}_2}\frac{2\alpha\beta^2\overline{\sigma}_2}{\gamma c}$; xii) $\frac{\beta}{\gamma}\alpha\beta c \gamma\underline{\sigma}_2\ge\frac{2\beta\overline{\sigma}_2}{\underline{\sigma}_2}\frac{4\alpha\beta^2\overline{\sigma}_2}{\gamma}$; xiii) $\frac{1}{2}\frac{\beta}{\gamma}\alpha\beta c\gamma\underline{\sigma}_2\ge\frac{1}{2}\frac{\beta}{\gamma}\alpha\beta c^2 \gamma\underline{\sigma}_2\ge\frac{2\beta\overline{\sigma}_2}{\underline{\sigma}_2}\frac{2\alpha\beta^2\overline{\sigma}_2 c}{\gamma}$; ixv) $\frac{1}{4}\alpha^2\beta^2\underline{\sigma}_2\gamma\underline{\sigma}_2\ge\frac{2\beta\overline{\sigma}_2}{\underline{\sigma}_2}\frac{\alpha^2\beta^2\overline{\sigma}_2^2}{c^2} $; xv) $\frac{1}{2}\alpha^2\beta^2\underline{\sigma}_2\gamma\underline{\sigma}_2\ge\frac{2\beta\overline{\sigma}_2}{\underline{\sigma}_2}\frac{2\alpha^2\beta^2\overline{\sigma}_2^2}{c}$; xvi) $\frac{1}{4}\frac{\beta}{\gamma}\alpha\beta c\gamma\underline{\sigma}_2\ge \frac{2\beta\overline{\sigma}_2}{\underline{\sigma}_2}\alpha^2\beta^2\overline{\sigma}_2^2$ by recalling that $c = 4\max\left\lbrace c_1,c_2,c_3,c_4 \right\rbrace\in\left[0,1\right] $ where $c_1 = \max\left\lbrace\alpha\overline{\rho},\alpha\gamma\overline{\sigma}_2 \right)\dot\max\left\lbrace 2, \frac{\overline{\sigma}_2}{\underline{\sigma}_2}\right\rbrace$, $c_2 = 2\beta \frac{\overline{\sigma}_2}{\underline{\sigma}_2}\max\left\lbrace2,\frac{\overline{\sigma}_2}{\underline{\sigma}_2}\right\rbrace$, $c_3 = \frac{\beta\overline{\sigma}_2}{\underline{\rho}}\max\left\lbrace2,\frac{2}{\alpha\gamma\underline{\sigma}_2},\frac{\overline{\sigma}_2}{\alpha\gamma\underline{\sigma}_2^2} \right\rbrace$ and $c_4 = 2\frac{\beta\overline{\sigma}_2^2}{\gamma\underline{\sigma}_2^2}\max\left\lbrace\frac{\overline{\sigma}_2}{\underline{\sigma}_2},\alpha\gamma\overline{\sigma}_2\right\rbrace$. The proof completes. $\hfill\blacksquare$

\section{ Proof of Lemma 4}

It is seen that $\mathcal{D} = -\Pi_{3}$ is with respect to the diagonal matrix $\varPsi = \text{diag}\left(\psi_1,\psi_2,\cdots,\psi_{p_2} \right)$ where $\psi_i\in\left[0,1\right]$, which leads $\mathcal{D}\left(\psi_1,\psi_2,\cdots,\psi_{p_2}\right)$ to be a convex combination of $2^{p_2}$ diagonal $p_2\times p_2$ matrices $\left\lbrace \mathcal{D}\left(k_1,k_2,\cdots,k_{p_2}\right) \mid k_i = 0 \text{ or } 1, i\in \mathbb{N}_{\left[1,p_2\right]}\right\rbrace$. Denote $\mathcal{D}_s$ the matrix with first $s$ elements in the diagonal being $1$ and the others $0$,where $s\in\mathbb{N}_{\left[0,p_2\right]}$. Without loss of generality, the lower bound of $\mathcal{D}_s$ implies the lower bound of $\mathcal{D}$. Note that
\begin{subequations}\nonumber
	\begin{align}
	\mathcal{D}_0 &= \left(\frac{\beta^2}{\gamma^2}-\frac{2\beta}{\gamma}\right)\alpha c\\
	&\ge 2\alpha^2\beta \left(1-\frac{\beta}{\gamma}\right)^2 A_2 A_2^\T \\
	&\ge \frac{1}{2}\alpha^2\beta A_2 A_2^\T\\
	\mathcal{D}_{p_2} &= \left(2-c \right)\alpha^2\beta A_2 A_2^\T \\
	&> \frac{1}{2}\alpha^2\beta A_2 A_2^\T
	\end{align}
\end{subequations}
In what follows, we write $A_2 A_2^\T$ into a block matrix form as follows
\begin{equation}\nonumber
A_2 A_2^\T = 
\begin{bmatrix}
\mathcal{A}_1& \mathcal{A}_2^\T\\
\mathcal{A}_2& \mathcal{A}_3
\end{bmatrix}
\end{equation}
where $\mathcal{A}_1\in\mathbb{R}^{s\times s}$, $\mathcal{A}_2\in\mathbb{R}^{\left(p_2-s\right)\times s}$ and $\mathcal{A}_3\in\mathbb{R}^{\left(p_2-s\right)\times \left(p_2-s\right)}$. Therefore, one can have
\begin{equation}
\begin{aligned}
\mathcal{D}_s &= 
\begin{bmatrix}
\alpha^2\beta\left(2-c\right)\mathcal{A}_1 & \alpha^2\beta\left(1-\frac{\beta}{\gamma}\right)\mathcal{A}_2^\T\\
\alpha^2\beta\left(1-\frac{\beta}{\gamma}\right)\mathcal{A}_2 & \alpha c\left(\frac{2\beta}{\gamma}-\frac{\beta^2}{\gamma^2}\right)\boldsymbol{I}
\end{bmatrix}\\
& \ge 
\begin{bmatrix}
\alpha^2\beta\left(2-c\right)\mathcal{A}_1 & \alpha^2\beta\left(1-\frac{\beta}{\gamma}\right)\mathcal{A}_2^\T\\
\alpha^2\beta\left(1-\frac{\beta}{\gamma}\right)\mathcal{A}_2 & 2\alpha^2\beta\left(1-\frac{\beta}{\gamma}\right)^2\mathcal{A}_3
\end{bmatrix}\\
& \ge \left(\frac{3}{2}-c\right)\alpha^2\beta A_2 A_2^\T
\end{aligned}
\end{equation}
where the first inequality follows from $c\left(\frac{2\beta}{\gamma}-\frac{\beta^2}{\gamma^2}\right)\boldsymbol{I}\ge 2\alpha\beta\left(1-\frac{\beta}{\gamma}\right)^2\overline{\sigma}_2\boldsymbol{I}\ge2\alpha\beta\left(1-\frac{\beta}{\gamma}\right)^2 \mathcal{A}_3$ according to (\ref{Pi_3_IECO_bound1}). The proof is completed. $\hfill\blacksquare$



\ifCLASSOPTIONcaptionsoff
  \newpage
\fi



%
\bibliography{mybibfile}

%








\end{document}